\documentclass{article}

\usepackage[english]{babel}

\usepackage{amssymb,amsthm}
\usepackage{amsxtra}
\usepackage{mathtools}
\usepackage{txfonts}
\mathtoolsset{showonlyrefs=true}
\usepackage{subcaption}
\usepackage{tikz-cd}
\usetikzlibrary{cd}
\usepackage{mathdots}

\usepackage[numbers]{natbib}
\usepackage{todonotes}

\usepackage{graphicx}
\usepackage[colorlinks=true, allcolors=blue]{hyperref}

\newcommand{\Ss}{\mathbb{S}}

\newcommand{\Rr}{\mathbb{R}}

\newcommand{\SU}{\mathfrak{su}}

\newcommand{\Tr}{\mbox{Tr}}

\newcommand{\OW}{\overline{W}}
\newcommand{\OP}{\overline{P}}
\newcommand{\TW}{\widetilde{W}}

\title{Sparse-stochastic model reduction for 2D Euler equations}
\author{Paolo Cifani\footnote{Gran Sasso Science Institute}\footnote{University of Twente - \url{p.cifani@utwente.nl}}, Sagy Ephrati\footnote{University of Twente - \url{s.r.ephrati@utwente.nl}} and Milo Viviani\footnote{Scuola Normale Superiore of Pisa - \url{milo.viviani@sns.it}}}
\date{}

\begin{document}
\maketitle

\begin{abstract}
The 2D Euler equations are a simple but rich set of non-linear PDEs that describe the evolution of an ideal inviscid fluid, for which one dimension is negligible.
Solving numerically these equations can be extremely demanding.
Several techniques to obtain fast and accurate simulations have been developed during the last decades.
In this paper, we present a novel approach which combines recent developments in the stochastic model reduction and conservative semi-discretization of the Euler equations.
In particular, starting from the Zeitlin model on the 2-sphere, we derive reduced dynamics for large scales and we close the equations either deterministically or with a suitable stochastic term.
Numerical experiments show that, after an initial turbulent regime, the influence of small scales to large scales is negligible, even though a non-zero transfer of energy among different modes is present.
\end{abstract}

\section{Introduction}
The 2D Euler equations are a fundamental model for ideal fluids \cite{Eu1757}.
During the last two centuries, these equations have stimulated an intense activity both in terms of mathematics and physics (see for example the seminal works of Helmholtz and Arnol'd \cite{Hel1858, arn}).
In computational science and numerical analysis, retaining at a discrete level the rich non trivial structure of these equations is still a challenging problem \cite{ara,MoVi2020}.
One main computational issue is the "curse of dimensionality".
Indeed, turbulent phenomena vary in different spatial and time scales and in order to capture them, high resolution simulations are needed.

A peculiar aspect of 2D ideal fluids is the presence of infinitely many conservation laws.
In particular, the conservation of energy and enstrophy (the $L^2$ norm of the curl of the velocity field) implies a double cascade phenomenon \cite{Kr1967}: the energy tends to move from small scales to large scales, whereas the enstrophy tends to follow the opposite direction.
Hence, in terms of the curl of velocity, or vorticity, it is possible to clearly separate two regimes: one slowly evolving at large scales and one fast at small scales.
Theoretically, the study of non-deterministic fluid models for different regimes have gained interest in the SPDE community \cite{FlPa2021}.
The equations studied in \cite{FlPa2021} and the results proved therein, show a precise connection between different space-time regimes with a reduced model for large scales.
Indeed, it is shown that a suitable model for large scales is given by the so called SALT equations \cite{Hol2015}, in which a transport noise term models the infinitesimal action of the small scales on the large ones.
Several numerical tests have shown the usefulness of the SALT equations as a powerful tool for model reduction \cite{CoCrHoPaSh2020, EpCiLuGe2023}.

However, defining precisely what large and small scales are is still an open problem.
In this paper, we present a criterion for defining large scales in terms of truncation of Fourier expansion.
We point out that other choices and interpretations of large and small scales are possible (see for example \cite{ModViv2022}).
Let us first introduce the governing equations for the vorticity field $\omega$, defined on the 2-sphere $\Ss^2$ embedded in $\Rr^3$:
\begin{equation}\label{eq:Euler_equations}
\begin{array}{ll}
&\dot{\omega} = \lbrace\psi,\omega\rbrace\\
&\Delta \psi = \omega.
\end{array}
\end{equation}
The Poisson bracket is defined as
\[ \lbrace\psi,\omega\rbrace:=\nabla \psi\cdot\nabla^\perp\omega\]
and the Laplacian is the Laplace--Beltrami operator on $\Ss^2$.
As mention above, equations \eqref{eq:Euler_equations} have infinitely many first integrals: energy $H(\omega)=\frac{1}{2}\int_{\Ss^2}\psi\omega$, Casimirs $C_n(\omega)=\int_{\Ss^2}\omega^n$, for $n\geq 1$ and angular momentum.
Understanding the role played by these invariants is still an open problem, especially for long-time evolution of the fluid \cite{DoDr2022}.

In order to gain numerical insight on this question, V. Zeitlin proposed a spatial discretiziation of \eqref{eq:Euler_equations}, which retains many of the first integrals above \cite{ze1,ze2}.
The Euler--Zeitlin equations are defined as follows:
\begin{equation}\label{eq:EZ_equations}
\begin{array}{ll}
&\dot{W} = [P,W]\\
&\Delta_N P = W.
\end{array}
\end{equation}
Here $W$ is a $N\times N$ skew-Hermitian matrix with zero trace, that is, an element of the Lie algebra $\SU(N)$.
The bracket $[P,W]$ is the usual matrix commutator and the discrete Laplacian $\Delta_N$ is defined such that its spectrum is a truncation of the spectrum of $\Delta$ \cite{hopyau}.
As mentioned above, the Euler--Zeitlin equations possess the following integral of motions: energy $H(W)=\frac{1}{2}\Tr(PW)$, Casimirs $C_n(W)=\Tr(W^n)$, for $n = 2,\ldots, N$ and angular momentum.
The core of the Zeitlin model is how the original vorticity $\omega$ and the discrete one $W$ are linked.
Indeed, the representation theory of $SU(2)$ provides a deep connection between the discrete Laplacian $\Delta_N$ and a particular basis $\lbrace T_{lm}\rbrace$ of $\SU(N)$, for $l=1,\ldots,N-1$ and $m=-l,\ldots,m$ \cite{hopyau,bhss}:
\begin{itemize}
    \item each $T_{lm}$ is an eigenvector of $\Delta_N$, with eigenvalue $-l(l+1)$,
    \item for each $N\geq 1$, there exists a linear map $p_N:C^\infty(\Ss^2)\rightarrow\SU(N)$,
    defined via the (real) spherical harmonics basis $\lbrace Y_{lm}\rbrace$ as $p_N(Y_{lm})=T_{lm}$, if and only if $l\leq N-1$,
    \item $\|p_N\lbrace\psi,\omega\rbrace - N^{3/2}[p_N\psi,p_N\omega]\|\rightarrow 0$, for $N\rightarrow\infty$, where the norm is the operator one.
\end{itemize}

The classical way to determine large and small scales is to choose a wave number $\overline{l}$ as a threshold for the large scales (see for example \cite{BoEc2012,CoCrHoPaSh2020}).
In this work, we propose the following criterion to set the threshold $\overline{l}$.
Consider a time scale in which the fluid's energy spectrum profile has reached a stationary state.
Then, typically (that is, out of equilibrium) the spectrum exhibits a double slope, which determines a kink at a certain wave number $\overline{l}$.
Then, we defined the large scales $\OW$ as the filtered vorticity with modes up to  $\overline{l}$, obtaining a banded matrix.
We propose three possible ways, both deterministic and stochastic, of closing the equations for $\OW$, by choosing different interaction with the small scales. 
Finally, we provide numerical tests to assess the different models introduced.

\section{Sparse-stochastic model reduction}
The Euler--Zeitlin equations \eqref{eq:EZ_equations} allow to study some typical features of the 2D fluids in the matrix language. 
In this section, we propose a way to reduce the complexity of the equations \eqref{eq:EZ_equations}, by defining from $W$ a sparse matrix $\OW$ which retains the relevant large scale information. 
Then, we show different ways of closing the equations for $\OW$, adding a suitable stochastic term.

In the Zeitlin model, the basis element $T_{lm}$ of $\SU(N)$ have non-zero entries only in the lower and upper $\pm m$ diagonal.
If we look at the anti-diagonals, instead, we are looking at the components determining the value of the vorticity field at certain latitude bandwidth on $\Ss^2$, as shown in Figure~\ref{fig:struct_W}.
\begin{figure}[hbt!]
\[W=\begin{pmatrix}
\textcolor{blue}{NORTH}, \textcolor{red}{m=0}  & \ldots & \textcolor{red}{\longrightarrow} & \ldots & \textcolor{red}{|m|=N-1}  \\
\vdots & \ddots & \ddots & \textcolor{red}{\nearrow} & \vdots\\
\textcolor{red}{\downarrow}  & \ddots &\textcolor{blue}{\searrow} & \ddots & \textcolor{red}{\downarrow}\\
\vdots & \textcolor{red}{\swarrow} & \ddots & \ddots & \vdots\\
\textcolor{red}{|m|=N-1} & \ldots & \textcolor{red}{\longrightarrow} & \ldots   & \textcolor{blue}{SOUTH}, \textcolor{red}{m=0} \\
\end{pmatrix}\]
\caption{Structure of the discrete vorticity $W$ in the Zeitlin model.}\label{fig:struct_W}
\end{figure}

The large scales are typically chosen to be the modes such that $l$ is smaller than a threshold level $\overline{l}$.
In the Euler--Zeitlin model, this corresponds to consider the banded matrices limited in the diagonals $\pm l\leq\overline{l}$ and then removing the components corresponding to $l>\overline{l}$.
The Poisson equation which defines the stream matrix $P$ preserves this sparsity structure, since the basis elements $T_{lm}$, the eigenvectors of the Laplacian, are themselves sparse.
However, the Lie bracket does not restrict to this space.
Indeed at each time-step we have to project the vector field into the right space.

Usually, we do not have any chance to guess the contribution of the small scales to the evolution of the large ones.
However, we expect that after an initial turbulent transition, the fluid exhibits two clearly separated spatial scales.
The hint for such a scenario is due to several numerical simulations of the Euler--Zeitlin equations \cite{BoEc2012, ModViv2022}.
Eventually, the energy profile reaches a fixed configuration with two slopes.
The first part of the spectrum represents the distribution of energy at large scales, whereas the second part the distribution of energy at small scales.
Typically, the separation  between large and small scales occurs at a wave number $\overline{l}\approx \sqrt{N}$.
For wave numbers lager than $\overline{l}$ the energy spectrum has the characteristic slope of $l^{-1}$, which is the one of white noise, see Figure~\ref{fig:initial_vorticty}.
\begin{figure}
\begin{tikzcd}[ column sep=tiny]
& \overline{W}:=\begin{pmatrix}
\overline{\omega}_{11}  & \ldots & \overline{\omega}_{1\overline{l}} & 0  & \dots  \\
\vdots  & \ddots & \ddots & \ddots  & \vdots\\
\overline{\omega}_{\overline{l}1}  & \ddots & \ddots & \ddots  & \overline{\omega}_{\overline{l}N}\\
0  &  \ddots & \ddots  &\ddots  & \vdots \\
\vdots  & \ddots  & \overline{\omega}_{N\overline{l}} & \ldots  & \overline{\omega}_{N N}   \\
\end{pmatrix}
&[1.5em] \\
W=\begin{pmatrix}
\omega_{11} & \ldots & \omega_{1N}\\
\vdots & \vdots & \vdots\\
\omega_{N1} & \ldots & \omega_{NN}\\
\end{pmatrix} \ar[ur, "\pi"] \ar[dr, "\pi^\perp"']  & \\
& \begin{pmatrix}
\widetilde{\omega}_{11} & \ldots & \widetilde{\omega}_{1N}\\
\vdots & \vdots & \vdots\\
\widetilde{\omega}_{N1} & \ldots & \widetilde{\omega}_{NN}\\
\end{pmatrix}\approx \sum_{l=\overline{l}+1}^{N-1} \sum_{m=-l}^l \beta^{lm}T_{lm}=:\widetilde{W}\\
\end{tikzcd}
\caption{Filtering of large scale and definition of random small scale vorticity $\widetilde{W}$, via the independent Brownian motions $\beta^{lm}$.}\label{fig:filer_large_scale}
\end{figure}
The universal nature of the small scales suggests a model reduction in terms of large scales evolution combined with a stochastic term contribution.
In Figure~\ref{fig:filer_large_scale}, we show the procedure to get the two new fields $\OW$ and $\TW$.
To define $\OW$, we introduce the orthogonal projection $\pi$ onto the modes $l\leq \overline{l}$.
The small scales field $\TW$ is defined as the linear combination of the basis elements $T_{lm}$, for $l>\overline{l}$ with coefficients $\beta^{lm}$ as independent Brownian motions, with mean and variance obtained from the high resolution DNS.
Application of the Kolmogorov--Smirnov and Anderson--Darling tests for normality to the high resolution data suggest that the distribution of the basis coefficients for $T_{lm}$, for $l > \bar{l}$, is Gaussian.

Hence, we define $\OW:=\pi W$ and $\TW:=\sum_{l=\overline{l}+1}^{N-1} \sum_{m=-l}^l \beta^{lm}T_{lm}$.
With these new fields, we essentially have three possible choices.
The first one consists of a deterministic closure simply via the projection of the vector field onto the large scales:
\begin{equation}\label{eq:EZ_equations_det}
\begin{array}{ll}
&\dot{\OW} =\pi [\OP,\OW]\\
&\Delta_N\OP = \OW.
\end{array}
\end{equation}
The second model is the enstrophy-preserving stochastic closure, which is up to the projection $\pi$ a type of SALT equation (see \cite{Hol2015}):
\begin{equation}\label{eq:EZ_equations_casimir_noise}
\begin{array}{ll}
&d\OW = \pi[\OP,\OW]dt +  \sum_{l=\overline{l}+1}^{N-1} \sum_{m=-l}^l \frac{1}{-l(l+1)}\pi[T_{lm},\OW]\circ d\beta^{lm}\\
&\Delta_N\OP = \OW.
\end{array}
\end{equation}
Finally, the third one is a energy-preserving stochastic closure:
\begin{equation}\label{eq:EZ_equations_energy_noise}
\begin{array}{ll}
&d\OW = \pi[\OP,\OW]dt +  \sum_{l=\overline{l}+1}^{N-1} \sum_{m=-l}^l \pi[\OP,T_{lm}]\circ d\beta^{lm}\\
&\Delta_N\OP = \OW.
\end{array}
\end{equation}
We recall that the symbol $\circ$ denotes the Stratonovich integral.
In the next section, we perform a numerical test for the three different models \eqref{eq:EZ_equations_det}, \eqref{eq:EZ_equations_casimir_noise}, \eqref{eq:EZ_equations_energy_noise}, comparing them with the high resolution DNS.

\section{Numerical simulations}
In this section, we show a numerical experiment to study the performance of the models proposed in the previous section.
The numerical experiment is conducted as follows.
We set the high resolution level at $N=128$.
Then we generate a random initial condition and we run a high resolution DNS.
 \begin{figure}[hbt!]
\centering
 \captionsetup{width=0.6\textwidth}
 \includegraphics[width=1\linewidth]{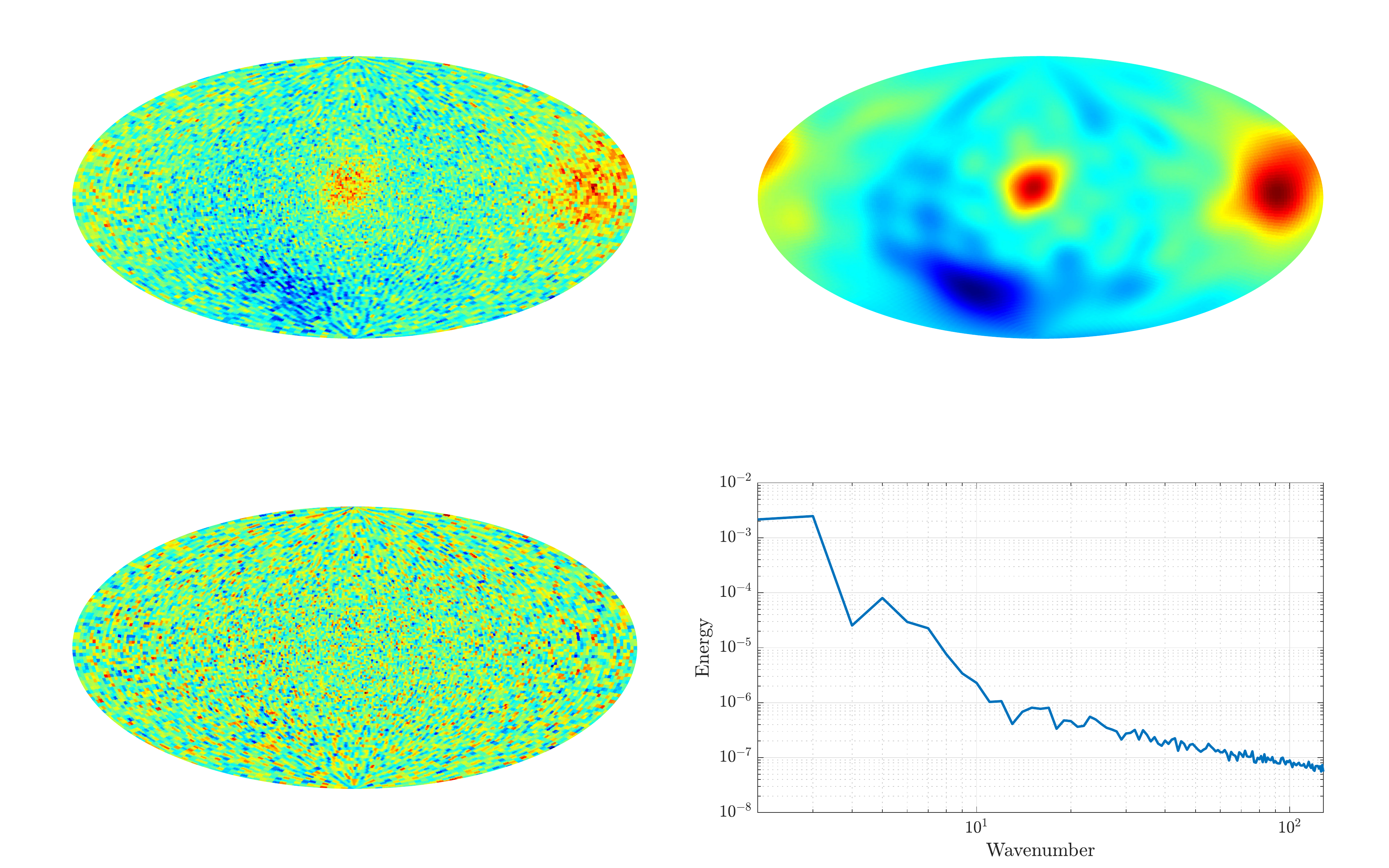} 
 \caption{Initial vorticity obtained via high resolution DNS. Top left, the field $W$, top right, the filtered field $\OW$, bottom left, $W-\OW$, bottom right, energy spectrum of $W$. Note the change of slope in the energy profile at $l\approx\sqrt{N}$.}
\label{fig:initial_vorticty}
 \end{figure}
We stop the simulation once a stationary energy profile is reached (see Figure \ref{fig:initial_vorticty}). 
Then, we select the large scale threshold as wave number $\overline{l}\approx\sqrt{N}$, at which the kink in the energy spectrum appears. 
In our numerical simulation the kink is found out to be at $\overline{l}=14$.
Then, we define our large scale field as $\overline{W}:=\pi W$, where $\pi$ denotes the orthogonal projection onto the modes for $l\leq\overline{l}$.
The projection consists of two steps: first we extract the components up to $\overline{l}$ and then we generate the field $\OW$.
The cost of calculating each component is $\mathcal{O}(N)$ and since we need to repeat this operation $\overline{l}^2-1\approx N$ times, the total cost of extracting the components is $\mathcal{O}(N^2)$.
Clearly, to construct the field $\OW$ we have to perform $\mathcal{O}(N^2)$ operations.
Hence, the total computational cost of the projection $\pi$ is $\mathcal{O}(N^2)$.
We also define $\TW$, as explained in the previous section.
Finally, we restart the original high resolution simulation and we perform numerical simulations where the small scales are modelled as described in equations \eqref{eq:EZ_equations_det}, \eqref{eq:EZ_equations_casimir_noise}, \eqref{eq:EZ_equations_energy_noise}, for approximately 250 time units.
In our numerical simulations, the time integration is done via the Heun-type scheme adapted for the SDEs, with time-step $h=0.25$.
\begin{figure}[hbt!]
\centering
 \captionsetup{width=0.6\textwidth}
 \includegraphics[width=1\linewidth]{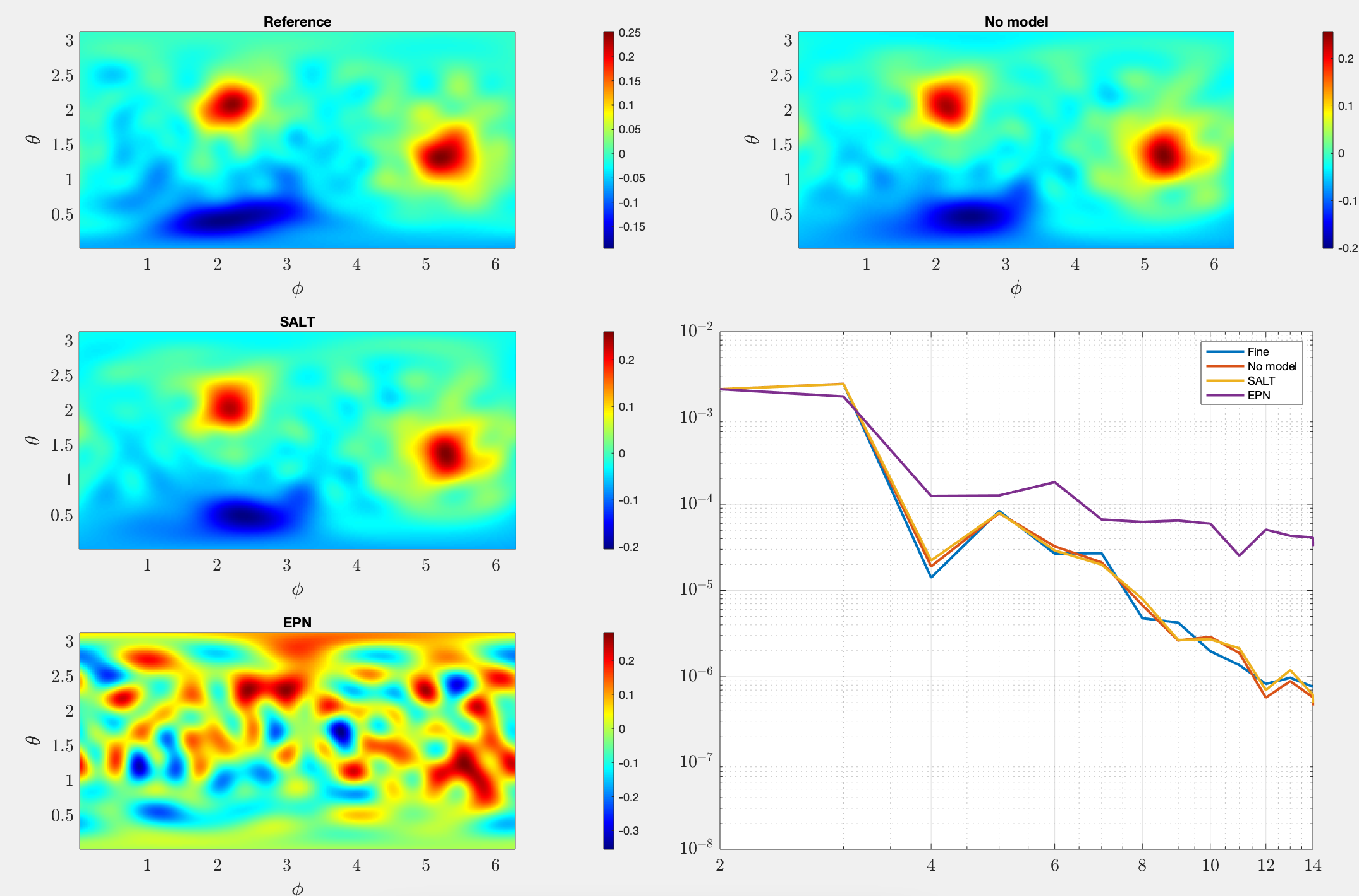} 
 \caption{Evolution and energy spectra of the large scales in 250 time units, via the different models proposed and the high resolution one. No model corresponds to \eqref{eq:EZ_equations_det}, SALT  to \eqref{eq:EZ_equations_casimir_noise} and EPN to \eqref{eq:EZ_equations_energy_noise}.}
\label{fig:simulations}
 \end{figure}
 
We notice from Figure~\ref{fig:simulations} that the no-model and SALT perform very well compared to the reference solution, both in spatial and energy spectrum profile.
On the contrary, the energy preserving scheme completely loses any accuracy and a cascade of energy to lower wave numbers occurs.
We can explain these facts in terms of energy flux among different modes.
In the energy preserving scheme, no energy can leave the large scales.
Hence, if the transfer of energy between different modes is non-zero, the conservation of the large scales energy prevents the energy to flow form large scales to small scales, causing an extra accumulation of energy $l\approx\overline{l}$.

\begin{figure}[hbt!]
\centering
 \begin{subfigure}{0.48\textwidth}
  \captionsetup{width=0.8\textwidth}
 \includegraphics[width=1\linewidth]{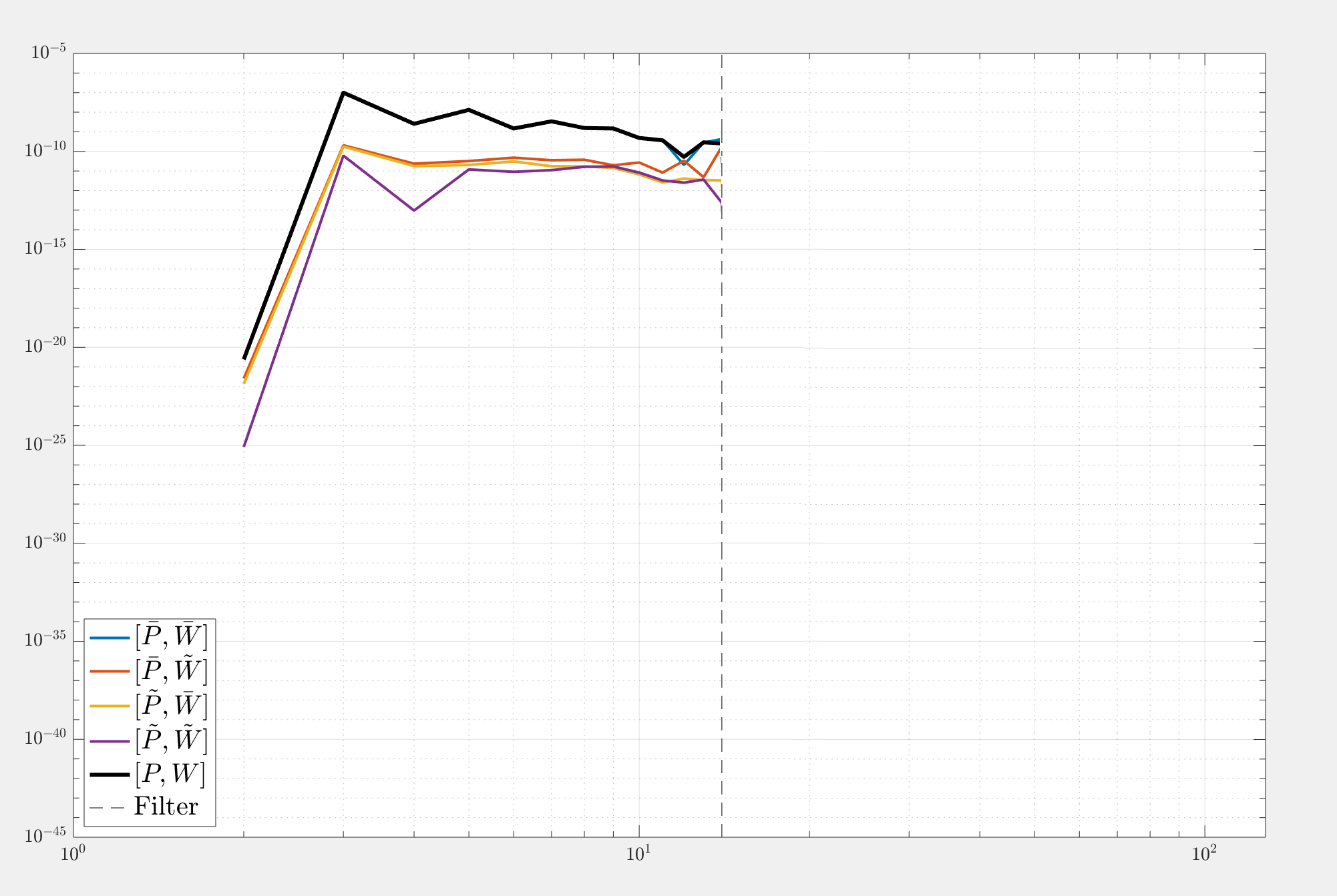} 
 \caption{Energy transfer among different modes at large scale, for high resolution DNS.}
\label{fig:EF_large_scales}
\end{subfigure}
 \begin{subfigure}{0.48\textwidth}
\captionsetup{width=0.8\textwidth}
 \includegraphics[width=1\linewidth]{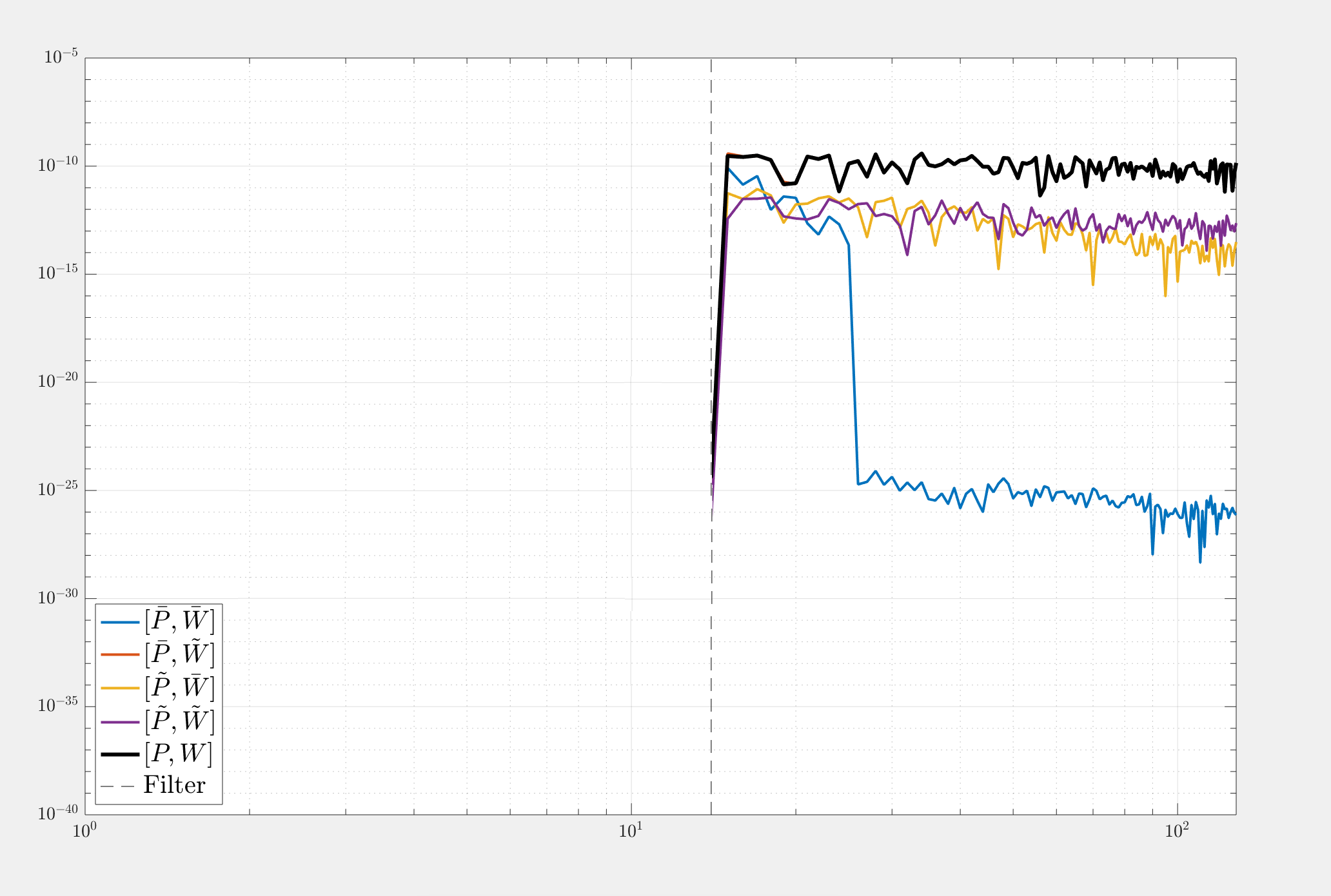} 
 \caption{Energy transfer among different modes at small scales, for high resolution DNS.}
\label{fig:EF_small_scales}
 \end{subfigure}
 \caption{}\label{fig:en_fluxes}
 \end{figure}

In order to check this thesis, we compute the energy transfer among different modes in the high resolution DNS.
Let us consider the energy at a level $l$:
\[E(l)=\frac{1}{2}\sum_{m=-l}^l \frac{\omega_{lm}^2}{l(l+1)}.\]
Then, the energy variation in time is given by 
\[\frac{dE(l)}{dt}=\sum_{m=-l}^l \frac{\omega_{lm}[P,W]_{lm}}{l(l+1)}.\]
Let $F(l):=|\frac{dE(l)}{dt}|$ be the absolute value of the energy transfer due to the non-linearity of the vector field $[P,W]$.
In Figure~\ref{fig:en_fluxes}, we plot the energy transfer contributions of the four possible coupling of large and small scales.
We notice that the transfer of energy between large and small scales is non-zero.
In particular, the main drivers of the energy for the components of $\OW$ is the vector field $[\OP,\OW]$, whereas for small scales is $[\OP,\TW]$.

\section{Conclusions and outlook}
In this paper, we have presented a possible strategy to reduce the complexity of the Euler--Zeitlin model, while performing long-time simulations. 
Numerical evidences show that the Euler--Zeitlin equations exhibit a clear separation of scales such that the large scale dynamics is quite robust to different coupling with small scales, either deterministic or stochastic.
Interestingly, the energy preserving scheme we have defined shows that the energy at large scales cannot be exactly conserved. This means that large and small scales are never completely decoupled, even when the energy spectrum profile reaches a stationary regime.
This indicate that for very long times a non-zero transfer of energy among different scales is present.

The Zeitlin model has been criticized for unrealistic conservation of enstrophy and other Casimirs at a finite level of truncation $N$.
Our result shows that this issue can be understood such that the Euler--Zeitlin equations are quite robust and precise in describing large scales, which means for wave numbers $l\approx\sqrt{N}$.
On the other hand, the remaining modes are themselves a model for the small scales, which correctly mimic the energy flux among different modes.

In conclusion, we have shown that the Zeitlin model can be a useful tool for simulating long-time large scale dynamics.
As future work, we aim to perform more systematic simulation using the parallelized code developed in \cite{CiMoVi2023} and available on \url{https://github.com/cifanip/GLIFS}.

\bibliographystyle{plainnat}
\bibliography{biblio.bib}
 
\end{document}